\documentclass[12pt]{article}
\usepackage{blindtext}
\usepackage{graphicx}
\usepackage[english]{babel}
\usepackage{amssymb}
\usepackage{float}
\usepackage{amsmath} 
\usepackage{mathrsfs}
\usepackage[ruled,vlined,linesnumbered]{algorithm2e}
\usepackage{algpseudocode}
\usepackage{multirow}
\usepackage{graphicx}
\usepackage{array}
\usepackage[left=2cm,right=2cm,
    top=2cm,bottom=2cm,bindingoffset=0cm]{geometry}
\usepackage{xcolor}
\definecolor{lightblue}{RGB}{30, 150, 186}
\usepackage[colorlinks,linkcolor=lightblue,urlcolor=lightblue,bookmarks=false,hypertexnames=true]{hyperref} 
\usepackage[singlelinecheck=false]{caption}
\usepackage{tocloft}

\begin{document}

\begin{center}
\textbf{AN IMPROVED KTNS ALGORITHM FOR THE JOB SEQUENCING AND TOOL SWITCHING PROBLEM}
\end{center}

\begin{center} Mikhail Cherniavskii, Department of Discrete Mathematics, Phystech School of Applied Mathematics and Informatics, Moscow Institute of Physics and Technology, Dolgoprudny, Russian Federation, e-mail: cherniavskii.miu@phystech.edu \end{center}
\begin{center} Boris Goldengorin, Department of Discrete Mathematics, Phystech School of Applied Mathematics and Informatics, Moscow Institute of Physics and Technology, Dolgoprudny, Russian Federation, e-mail: goldengorin.bi@mipt.ru \end{center}

\begin{center}
\textbf{ABSTRACT}
\end{center}
We outline a new Max Pipe Construction Algorithm (MPCA) with the purpose to reduce the CPU time for the classic Keep Tool Needed Soonest (KTNS) algorithm. The KTNS algorithm is applied to compute the objective function value for the given sequence of jobs in all exact and approximating algorithms for solving the Job Sequencing and Tool Switching Problem (SSP). Our MPCA outperforms the KTNS algorithm by at least an order of magnitude in terms of CPU times. Since all exact and heuristic algorithms for solving the SSP spend most of their CPU time on applying the KTNS algorithm we show that our MPCA solves the entire SSP on average 59 times faster for benchmark instances of D compared to current state of the art heuristics.

\begin{center}
\textbf{INTRODUCTION}
\end{center}

Currently as far as we are aware a globally optimal solution to the Job Sequencing and Tool Switching Problem (SSP) Tang, Denardo 
{[\ref {TangDenardo}]} is limited to a couple of dozen jobs with no more than 25 tools in each of them da Silva et al. {[\ref {daSilva2021}]}. This class of problems includes many optimization versions of SSP solved for a large number of small-scale orders (jobs) on conveyor lines (flexible manufacturing system - FMS). We do not consider technical and technological changes (improvements) in the equipment of the FMS in order to increase its productivity. Our goal is to increase the productivity of the conveyor line by finding an optimal sequence of loading the necessary tools, materials and other production resources (hereinafter referred to as tools) including the professional skills of employees, e.g. in order to minimize the conveyor's downtime. In other words, the methods and means of increasing the productivity of the conveyor in this article are related to an optimal sequencing of jobs each of which requires its specific collection of tools. Since a production line has a limited number of slots (magazine with a limited capacity) a tool switch is necessary to store the tools required to complete all jobs. Here by {\it tool switch} we understand the removal of a tool from the magazine (collection of slots) and the insertion of another tool in its place. Thus, the SSP consists in finding a sequence of jobs minimizing the total number of tool switches.

In this paper we assume that the following refined assumptions are satisfied Calmels {[\ref{Calmels}]}, da Silva et al. {[\ref {daSilva2021}]}:

$\bullet$ There is a set of jobs to be processed and each job requires a fixed set of specific tools; 

$\bullet$ The set of jobs and the subset of tools required for each job is known in advance;

$\bullet$  No job requires a set of tools that exceeds the capacity of the machine's magazine; 

$\bullet$  All tools are always available at least outside of the magazine;

$\bullet$  A single machine is available to process all jobs;

$\bullet$  It is an offline version of the problem;

$\bullet$  Once the machine has started processing a job it must be completed;

$\bullet$ Exactly one job can be processed at any time unit;

$\bullet$ The processing and completion times of a job are not dependent on the set of tools and do not impact the number of tool switching times;

$\bullet$ The tool sockets (slots) are identical;

$\bullet$  Only one tool switch is done at time;

$\bullet$  Each tool fits in any slot of the magazine and occupies one slot;

$\bullet$  The time associated with removal and insertion (switch) of a tool is independent and constant;

$\bullet$  No breaks, no wear and no maintenance of the tools are considered.

Here, a single machine is an abstraction. It reflects many different interpretations, e.g. huge and small companies, production and service lines, and organizations - educational, medical, governmental and private. For each of these we are able to indicate an input of jobs, their operations including available resources (tools) and an expected output of ordered (sequenced) products (jobs) to be processed Goldengorin, Romanuke {[\ref {GoldRom2021}]}.

In this paper we are not going to overview the state of the art for mathematical models, methods and algorithms for modeling and solving the SSP. We rather refer the interested reader to the recent literature review Calmels {[\ref{Calmels}]}, as well as the article da Silva et al. {[\ref {daSilva2021}]}. We would like to emphasize that the SSP is one of the NP-hard problems Tang, Denardo {[\ref {TangDenardo}]}. The relevant mathematical models and solution methods for SSP should take into account that the number of tool switches for the next job depends not only on a single or pair of prior scheduled jobs but, in the worst case, on all jobs scheduled before the pending job Crama et al. {[\ref{Crama}]}, Ghiani et al. {[\ref{Ghiani}]}, Ahmadi et al. {[\ref{AhmadiGoldengorin}]}. Most publications report globally optimal solutions to the entire SSP. However, their conclusions are questionable without any proof of how they take into account the real number of switches. For details we refer to Tables 6 and 7 in da Silva et al. {[\ref {daSilva2021}]}.

The SSP can be formulated as follows. We are given the sets of jobs $ J = \{1, \dots, n \} $ and tools $ T = \{1, \dots, m \} $, a single magazine $ C $ denoting the maximum number of tools (slots) that can be placed (occupied) in the machine's magazine, $C < m$, the set of tools  $ T_i \subset T $, $ i = 1, \dots n  $ required by job $i$ which should be in the magazine in order to complete the job $i$, $ | T_i | \leq C $. Thus, for the given sequence of $n$ jobs we are going to associate $ n $ instants each of which will be presented by the set tools $ M_i \subset T $ in the magazine. A feasible SSP solution is the sequence $ M = \{M_i \} _ {i = 1} ^ {n} $, such that sequence of jobs $(1, 2, \dots , n)$ can be completed. An optimal solution to the SSP is a  sequence $ M $ that minimizes the total number of tool switches required to move from one job to another and complete all jobs. 

As shown by Tang, Denardo {[\ref{TangDenardo}]} SSP can be decomposed into 2 following problems.  
\begin{enumerate}
\item Tool Loading Problem(TLP) - for a given sequence of jobs, find the optimal sequence of magazine states $M$ that minimizes the total number of tool switches.
\item Job Sequencing Problem(JeSP) – finding a sequence of jobs such that the number of tool switches is minimal after solving the TLP for this sequence. 
\end{enumerate}

Note that most publications devoted to solving the SSP consider heuristics based on completely different classes of metaheuristics, e.g., tabu search, iterated local search, and genetic algorithms {[\ref {Mecler}]}. Regardless of the designed heuristic's nature they try to replace all permutations defined on the entire set of jobs and select the best values of the SSP objective function (OF), i.e. the minimum number of tool switches. One of the first and most popular algorithms to compute the SSP OF, well known for over thirty years, is the Keep Tool Needed Soonest (KTNS) algorithm {[\ref{TangDenardo}]}. The main purpose of our paper is to improve the KTNS algorithm since the efficiency of any exact or heuristic algorithm is based on a partial enumeration of SSP feasible solutions and depends on the CPU time to compute the SSP OF value.

Our paper is organized as follows. Before describing our new Max Pipe Construction Algorithm (MPCA) we provide a small numerical example to illustrate the pitfalls of slowing down the SSP OF value computation. In the next two sections we describe the solution method and formulate statements to justify the correctness of our MPCA. We provide an example demonstrating the inner workings of MPCA and design an efficient $ MPCA-bitwise $ implementation and evaluate its time complexity. The Experimental Results section includes our computational study and the final section contains a summary and future research directions.

\begin{center}
\textbf{A Numerical Example}
\end{center}

\begin{figure}[H]
\begin{minipage}[t]{.45\linewidth}
  \includegraphics[width=\linewidth]{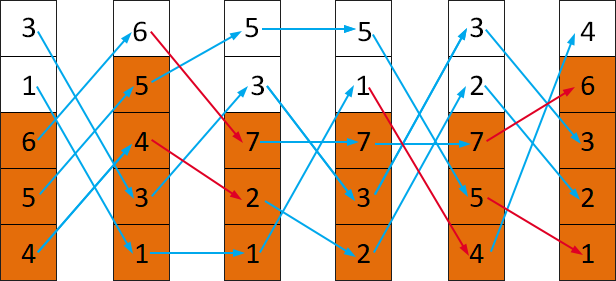}\label{figure:solution_example1}
  \caption{Solution of the \textit {TLP} problem for the $ (1,2,3,4,5,6) $ sequence.}
\end{minipage}
\begin{minipage}[t]{.1\linewidth}
$ $
\end{minipage}
\begin{minipage}[t]{.45\linewidth}
  \includegraphics[width=\linewidth]{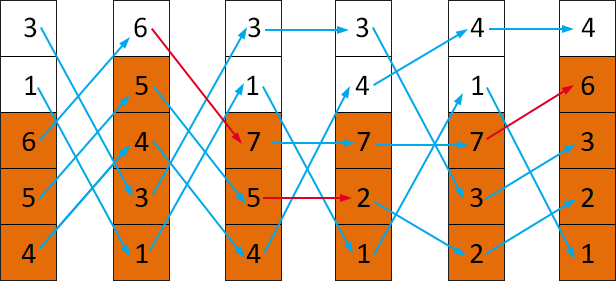}\label{figure:solution_example2}
  \caption{Solution of the \textit{TLP} problem for the $ (1, 2, 5, 3, 4, 6) $ sequence.}
\end{minipage}
\end{figure}

Let's consider an example. $T_1=\{4,5,6\},T_2=\{1,3,4,5\},T_3=\{1,2,7\},T_4=\{2,3,7\}
,T_5=\{4,5,7\},T_6=\{1,2,3,6\}$ is given sequence of jobs represented as sets of tools that they require, $C=5$ is magazine capacity. \hyperref[figure:solution_example1]{Fig. 1} shows the optimal loading of empty slots to perform job in the order of $ (1,2,3,4,5,6) $. \hyperref[figure:solution_example2]{Fig. 2} shows the optimal loading of empty slots for performing job in the order $ (1, 2, 5, 3, 4, 6) $. The sequence $ (1, 2, 5, 3, 4, 6) $ is found by exhaustive brute force search and is a solution to the problem \textit {JeSP}. In the \hyperref[figure:solution_example2]{Fig. 1,2}, red arcs represent tool switches. Thus, the smallest number of switches for the $ (1,2,3,4,5,6) $ sequence is $ 5 $, and for the $ (1, 2, 5, 3, 4, 6) $ sequence is $ 3 $.

\noindent Crama et al. [\ref{Crama}] proved that the \textit{JeSP} is NP-hard. Usually, TLP is solved with purpose to calculate the SSP objective function value  for each job sequence in \textit{JeSP}. To solve TLP Tang, Denardo {[\ref {TangDenardo}]} proposed the \textit{Keep Tool Needed Soonest (KTNS)} with time complexity $ O(mn) $. Ghiani et al. {[\Ref{Ghiani}]} modified \textit{KTNS} to \textit{Tailored KTNS} with the purpose to calculate a lower bound to the number of switches between two jobs which was used to accurately solve \textit{JeSP}.

\begin{center}
\textbf{PROBLEM FORMULATION}
\end{center}

Let's refresh the basic notation for  \textit{the job Sequencing and tool Switching Problem (SSP)} and further illustrate them with examples.
	$T = \{ 1,2, \dots ,m \}$ is the set of tools, where $ m $ is the number of tools.
	$J=\{1,2 \dots n \}$ is the set of jobs, where $ n $ is the number of jobs.
	$T_i \subset T$ is the set of tools needed to do the job $ i $, where $ i \in J $. $ C <m $ - capacity (capacity) of the magazine.
	$M_i \subset T$ is the state of the magazine when the job $ i $ is performed i.e. the set of tools located in the magazine at the time of the job $ i $, where $i \in J$, $|M_i|=C$, $T_i \subseteq M_i$.
	$\boldsymbol{S} = \{ S: S = (\sigma(1), \sigma(2), \dots , \sigma(n) )$, where $ \sigma $ is  a permutation $ \} $ is the set of all reordering of jobs. 
	$\boldsymbol{M}(S) = \{ M = (M_{1}, \dots M_{{n}}): T_i \subseteq M_i, |M_i|=C, i \in J \}$  is the set of all sequences of magazine states such that jobs can be performed in order $S \in \boldsymbol{S}$. 
	$\boldsymbol{M} = \{ M \in \boldsymbol{M}(S): S \in \boldsymbol{S} \}$ is the set of all possible magazine states.
	$switches(M) = \sum_{i=1}^{n-1} C - |M_{{i}} \cap M_{{i+1}}|$ is the number of switches for the sequence of magazine states, where $M \in \boldsymbol{M}$. 
The $ TLP $ problem is formulated as finding $\underset{M \in \boldsymbol{M}(S)}{argmin} \{ switches(M) \}$. The problem $ JeSP $ is formulated as finding $\underset{S \in \boldsymbol{S}}{argmin} \{
	\underset{M \in \boldsymbol{M}(S)}{min} \{ switches(M) \}
\}$.

\begin{table}[H]

  \caption{
  Example of $KTNS$ processing.
  }

  \begin{tabular}{ | c | m{9cm} | }

    \hline
    \begin{minipage}{0.5\textwidth}
      \includegraphics[width=\linewidth]{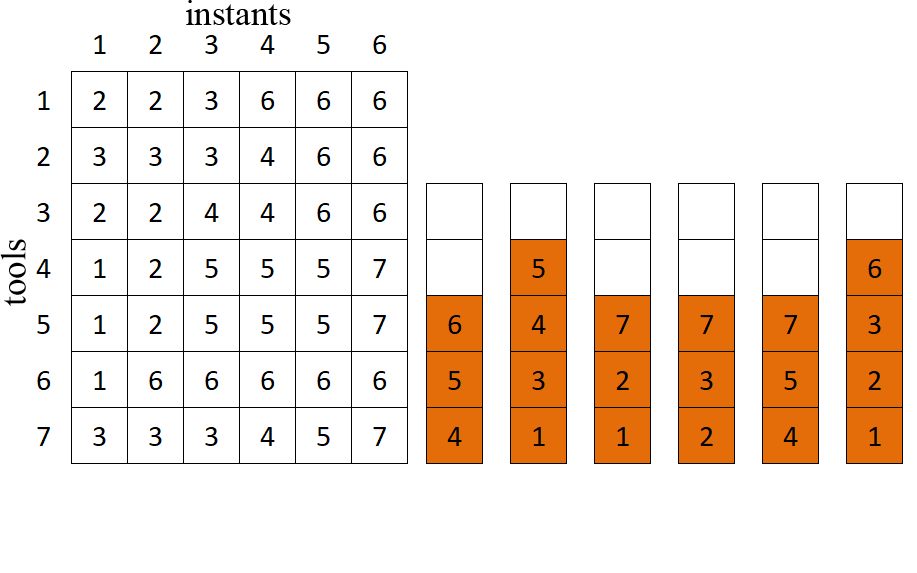}
    \end{minipage} &
$T_1=\{4,5,6\},T_2=\{1,3,4,5\},T_3=\{1,2,7\},T_4=\{2,3,7\}
,T_5=\{4,5,7\},T_6=\{4,5,7\}$ is given sequence of jobs represented as sets of tools that they require(marked in orange).
At first KTNS calculates an auxiliary matrix, the element $ a_ {ij} $ of which denotes the first instant starting from $ j $ at which the tool $ i $ is needed. For example starting from instant $2$ the first instant at which the tool $5$ is needed is instant $2$ itself, the first instant in which the tool $6$ is needed is $6$. Starting from instant $6$ the tool $4$ is never needed again than $a_{4,6} = 7$ 
    \\
    
    \hline
    \begin{minipage}{0.5\textwidth}
      \includegraphics[width=\linewidth]{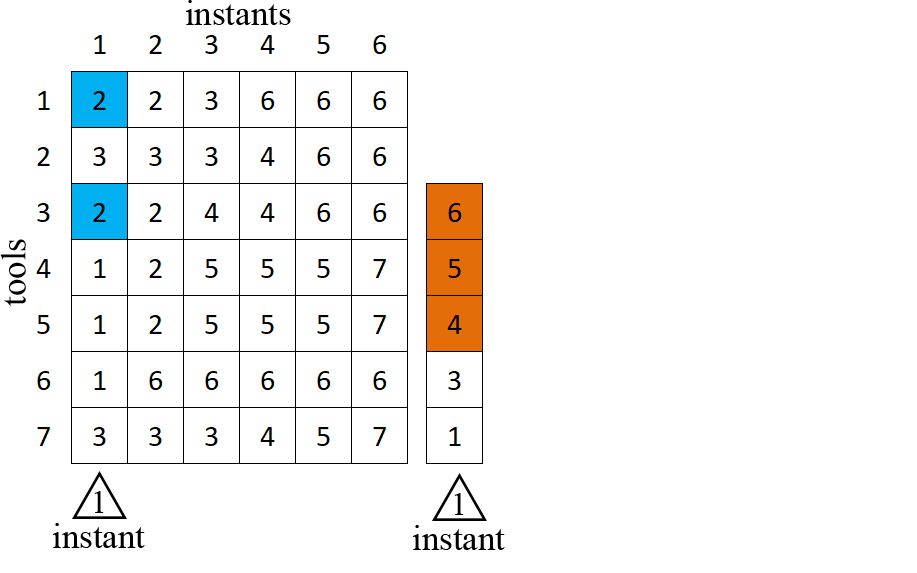}
    \end{minipage} &
Let's consider the first point in time. Tools required for the job scheduled at this moment $ T_1 = \{4,5,6 \} $ (marked in orange). We fill 3 slots with these tools. The remaining 2 slots are filled with the tools that are needed soonest. We see in the table (column 1) that the $ 1 $ tool will be needed at the moment $ 2 $, the $ 2 $ tool will be needed at the moment $ 3 $, the $ 3 $ tool will be needed at the moment $ 2 $, the $ 7 $ tool will be needed at the moment $ 3 $. Choose the $ 1 $ and $ 3 $ tools since they are needed soonest. 
    \\
	\hline
	
    \begin{minipage}{0.5\textwidth}
      \includegraphics[width=\linewidth]{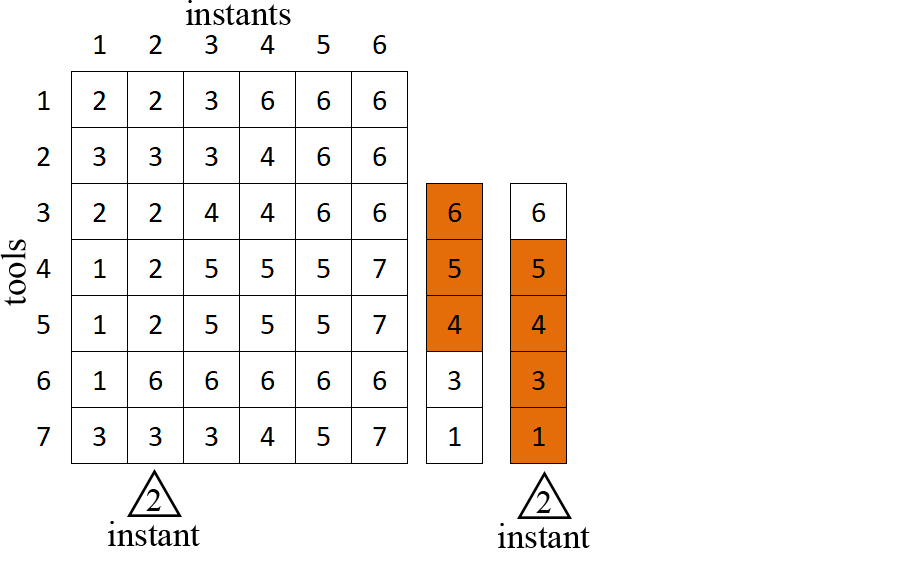}
    \end{minipage} &
Let's consider the instant $ 2 $. The tools necessary for the job $ T_2 = \{1,3,4,5 \} $ scheduled at this moment (marked in orange), all these tools were already present at the previous step, so no switches are required. Let's move on to the next step.
    \\
	\hline

  \end{tabular}
\label{table:MPCA_example_2}
\end{table}

\begin{table}[H]

  \caption{
   Example of $KTNS$ execution.
  }

  \begin{tabular}{ | c | m{9cm} | }
  
  \hline
      \begin{minipage}{0.5\textwidth}
      \includegraphics[width=\linewidth]{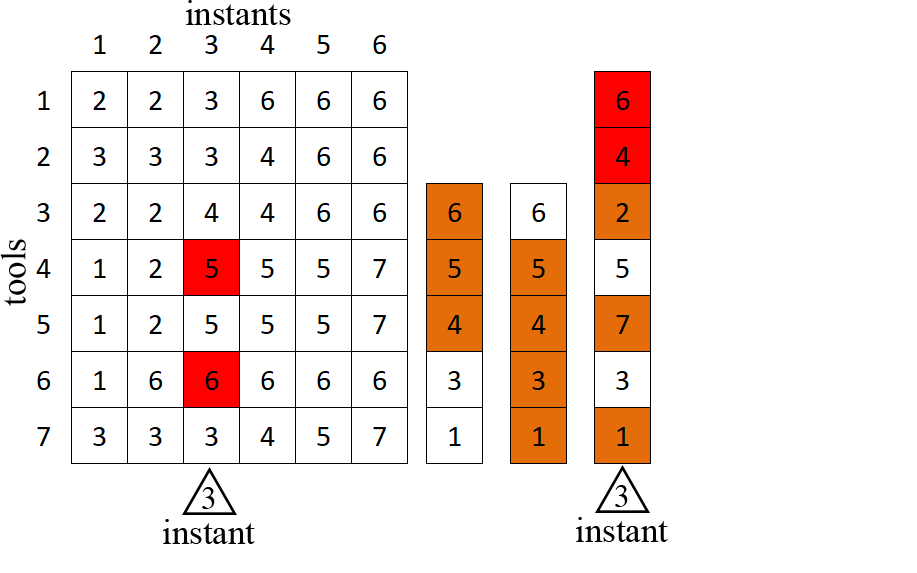}
    \end{minipage} &
Let's consider the instant $ 3 $. Tools required for the job scheduled at this moment $ T_3 = \{1,2,7 \} $ (marked in orange). The $ 1 $ tool was already in the magazine in the previous step, and the $ 7 $ and $ 2 $ tools need to be charged, but the magazine is already full, so you need to remove two tools. Let's remove the tools that are needed later. We see in the table (column 3) that the $ 3 $ tool will be needed at the moment $ 4 $, the $ 4 $ tool will be needed at the moment $ 5 $, the $ 5 $ tool will be needed at the moment $ 5 $, the $ 6 $ tool will be needed at the moment $ 6 $. Let's remove the $ 4 $ and $ 6 $ tools as they will be needed later than others. Thus, there were two switches of the tool $ 6 $ for the tool $ 2 $ and the tool $ 4 $  for the tool $ 7 $.
    \\

\hline
    \begin{minipage}{0.5\textwidth}
      \includegraphics[width=\linewidth]{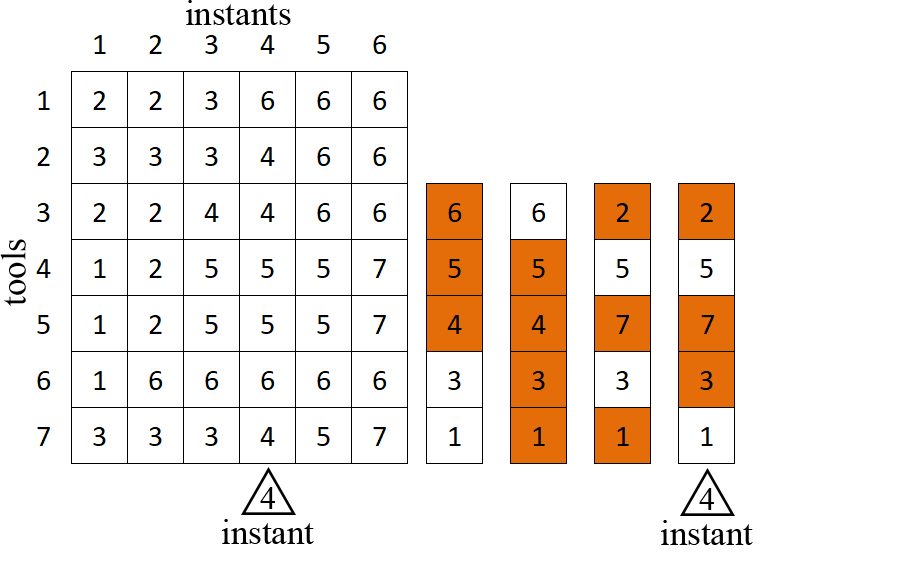}
    \end{minipage} &
Let's consider the instant $ 4 $. The tools required for the job $ T_4 = \{2,3,7 \} $ scheduled at this moment (marked in orange), all these tools were already present at the previous step, so no switches are required. Let's move on to the next step.
    \\
	\hline

    \begin{minipage}{0.5\textwidth}
      \includegraphics[width=\linewidth]{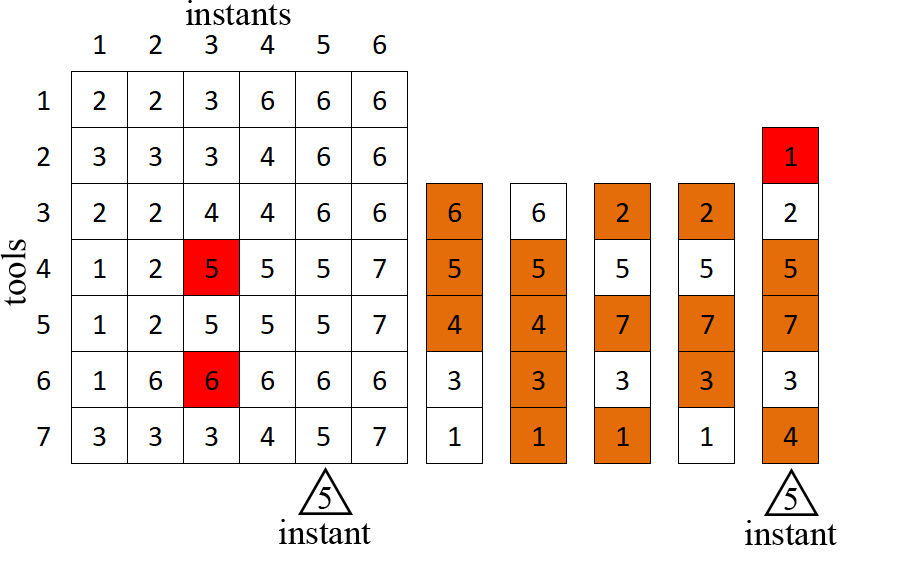}
    \end{minipage} &
Let's consider the instant $ 5 $. The tools required for the job scheduled at this moment $ T_5 = \{4,5,7 \} $ (marked in orange), The $ 5,7 $ tools were already in the magazine in the previous step, and the tool $ 4 $ needs to be charged, but the magazine is already full , so we need to delete one tool. Let's delete the tool that will be needed later than others. We see in the table (column 4) that the $ 1 $ tool will be needed at the instant $ 6 $, the $ 2 $ tool will be needed at the instant $ 6 $, the tool $ 3 $ will be needed at the moment $ 6 $. Let's remove the tool $ 1 $ as it is needed later then others. Instead of $ 1 $, it would be possible to delete $ 2 $ or $ 3 $, which would not affect the total number of switches. Thus, the $ 1 $ tool was switched with a $ 4 $ tool.
    \\
	\hline

  \end{tabular}
\label{table:KTNS_example_2}
\end{table}

\begin{table}[H]

  \caption{
   Example of $KTNS$ execution.
  }

  \begin{tabular}{ | c | m{9cm} | }
	\hline

    \begin{minipage}{0.5\textwidth}
      \includegraphics[width=\linewidth]{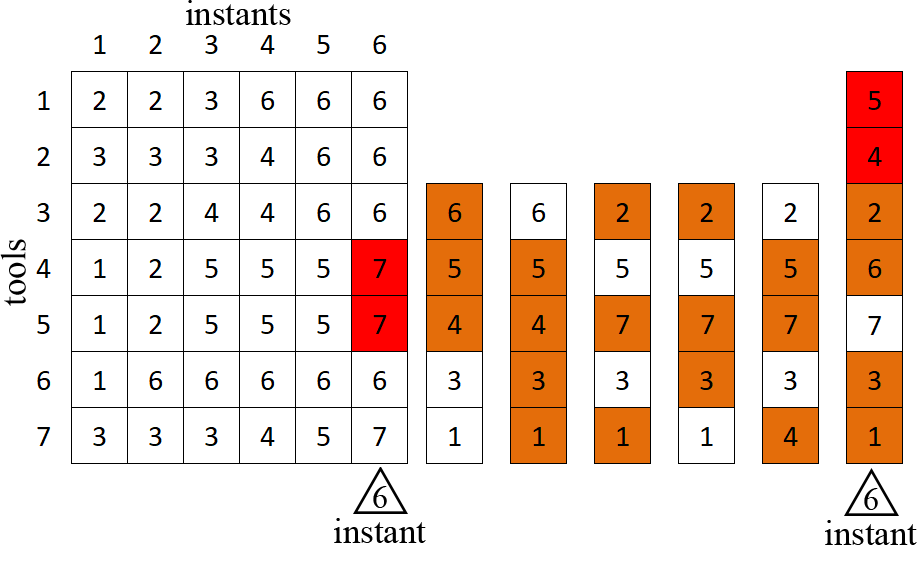}
    \end{minipage} &
Let's consider the instant $ 6 $. The tools needed for the job scheduled at this moment $ T_6 = \{1,2,3,6 \} $ (marked in orange), The $ 2,3 $ tools were already in the magazine in the previous step, and the $ 1,6 $ tools need to be charged, but the magazine is already full, so two tools must be removed. Let's remove the tools that are needed later than others. We see in the table (column 6) that the $ 5 $ tool will be needed at the time $ 7 $, the $ 4 $ tool will be needed at the time $ 7 $, the $ 7 $ tool will be needed at the time $ 7 $. Let's remove the $ 4.5 $ tools as they will be needed later then others. Instead of removing $4,5$, it would be possible to remove any pair from $4,5,7$ without affecting the total number of switches. Thus, there were two replacements of the tool $ 4 $  for the tool $ 6 $  and the tool $ 4 $  for the tool $ 1 $. KTNS has finished job with $5$ switches.
    \\
	\hline

  \end{tabular}
\label{table:KTNS_example_2}
\end{table}


\begin{center}
\textbf{SOLUTION METHOD}
\end{center}

\begin{figure}[H]
\begin{center}
  \includegraphics[width=.5\linewidth]{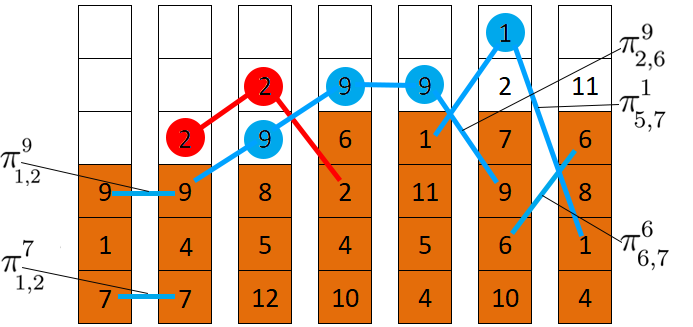}\label{figure:solution_example2}
  \caption{Examples of pipes.}
\end{center}
\end{figure}

Let's define $\mathscr{T}(M) = \{\pi_{s,e}^{tool}: tool \in (T_{{s}} \cap T_{{e}}) \backslash (\bigcup_{i=s+1}^{e-1} T_{i}), tool \in  \bigcap_{i=s+1}^{e-1} M_{i} \}$ as the set of all such triplets $(s,e,tool)$, that the $tool$ is used for job at moments $ s, e $, not used for job at moments $s+1, \dots , e-1$, however, is present in the magazine at moments $s+1, \dots , e-1$, despite the fact that at these moments the $ tool $ is not used for job. Let's Consider the main objects of this job - the elements of the set $\mathscr{T}(M)$, we will call them $pipes$, denote $\pi_{start,end}^{tool}$ - pipe with start at $ start $, end at $ end $ and tool $ tool $. \hyperref [figure: solution_example1] {Fig. 3} shows examples of pipes. Informally, the pipe is the saving of the $ tool $  from the moment $ start $, where it was used for job (marked in orange) until the moment $ end $, where it will again be used for job (marked in orange), but at intermediate points in time $ start + 1, \dots end-1 $ it must not be used to job (marked in white). So for example, $ \pi_ {5,7} ^ {1} $ saves tool $ 1 $ from time $ 5 $ to time $ 7 $ even though it is not used for job at time $ 6 $. Note that $ \pi_{1,2}^{7} $ is also a pipe, although there are no intermediate times between the times $ 1,2 $. The figure marked with a red pipe is not, since at the moment of time $ 1 $ the tool $ 2 $ is not used for job, i.e. $ 2 \notin T_{{1}} $. Let's call the capacity of the pipe the number of empty slots that is necessary for the existence of the pipe, which is equal to $end-start-1$.

\begin{algorithm}[H]
\SetAlgoLined
\DontPrintSemicolon
\SetKw{KwGoTo}{go to}
	$pipes\_count:= 0$\;	
	$M_1: = T_1, M_2: = T_2, \dots , M_n: = T_n$\\
	\For{$end=2,\dots,n$}{
		\For{$start=end-1 , \dots ,1 $}{
			$candidates:=\{\pi_{start,end}^{tool}: tool \in (T_{{start}} \cap T_{{end}}) \backslash (\bigcup_{i=s+1}^{e-1} T_{i}) \}$\;
			$empty\_slots:=min\{C- |M_i| :i \in \{start+ , \dots , end-1\} \}$\;
			\If{$|candidates| > empty\_slots$}{
				$candidates:=$ arbitrary elements from $candidates$ in the amount of $empty\_slots$.
			}
			\For{$\pi_{start,end}^{tool} \in candidates$}{
				add $tool$ to $M_{start+1} , \dots , M_{start+2} , \dots , M_{end - 1}$\;
			}
			$pipes\_count:= pipes\_count + |candidates|$\;
		} 
    }
    \Return $pipes\_count$ 
\

  \caption{\textit{MPCA} }
\end{algorithm}\label{algorithm:GPCA1}

$ $ \\
$ $ \\

\noindent  \textbf{Theorem 1.} 
Let $ C $ be the capacity of the magazine, $ T_1, \dots T_n $ are the required sets of tools for jobs $ 1, \dots, n $, $ S \in \boldsymbol {S} $ is the sequence of jobs, then
$$\underset{M \in \boldsymbol{M}(S)}{\mathrm{min}}\{switches(M)\} =  - \underset{M \in \boldsymbol{M}(S)}{\mathrm{max}}\{|\mathscr{T}(M)|\} -C + \sum_{i=1}^{n}|T_i| .$$

\noindent  \textbf{Theorem 2.} 
Let $ C $ be the capacity of the magazine, $ T_1, \dots T_n $ are the required sets of tools for jobs $ 1, \dots, n $, $ S \in \boldsymbol {S} $ is the sequence of jobs, $S \in \boldsymbol{S}$ - sequence of jobs, then
$$MPCA(S) = MPCA(T_1, \dots, T_n; \ C; \ S ) = \underset{M \in \boldsymbol{M}(S)}{\mathrm{max}}\{|\mathscr{T}(M)|\}.$$

\noindent From \hyperref[theorem:th1]{Theorem 1}, \hyperref[theorem:th2]{Theorem 2} it follows that the value of the objective function in the problem \textit{JeSP} equals  $\underset{M \in \boldsymbol{M}(S)}{\mathrm{min}}\{switches(M)\} =  - MPCA(S) -C + \sum_{i=1}^{n}|T_i|$.

\begin{table}[H]

  \caption{
   Example of $MPCA$ execution.
  }

  \begin{tabular}{ | c | m{9cm} | }
    \hline
    
    \begin{minipage}{0.4\textwidth}
      \includegraphics[width=\linewidth]{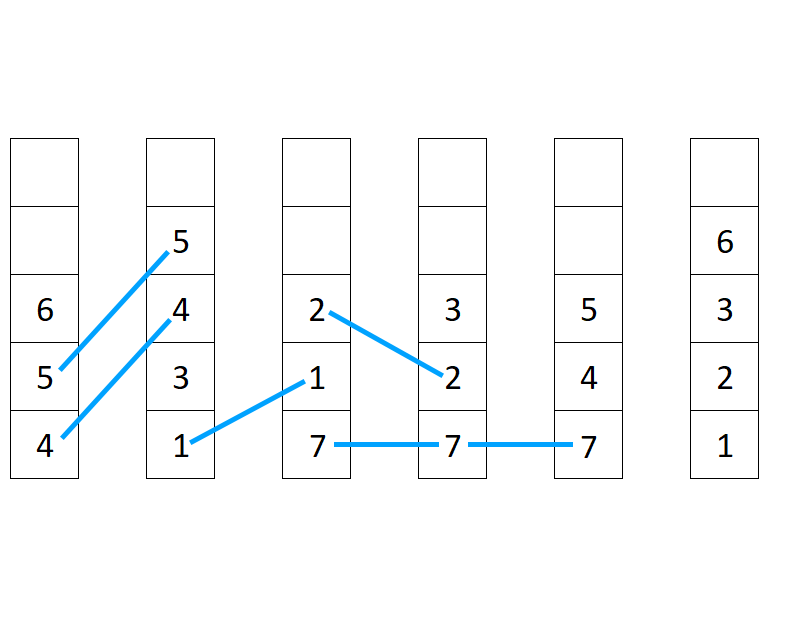}
    \end{minipage} &
$T_1=\{4,5,6\},T_2=\{1,3,4,5\},T_3=\{1,2,7\},T_4=\{2,3,7\}
,T_5=\{4,5,7\},T_6=\{1,2,3,6\}$ is given sequence of jobs represented as sets of tools that they require.
    Let's depict all pipes of capacity $0$. What is implemented in the algorithm line $ 1 $.
    
    \\
	\hline

    \begin{minipage}{0.4\textwidth}
      \includegraphics[width=\linewidth]{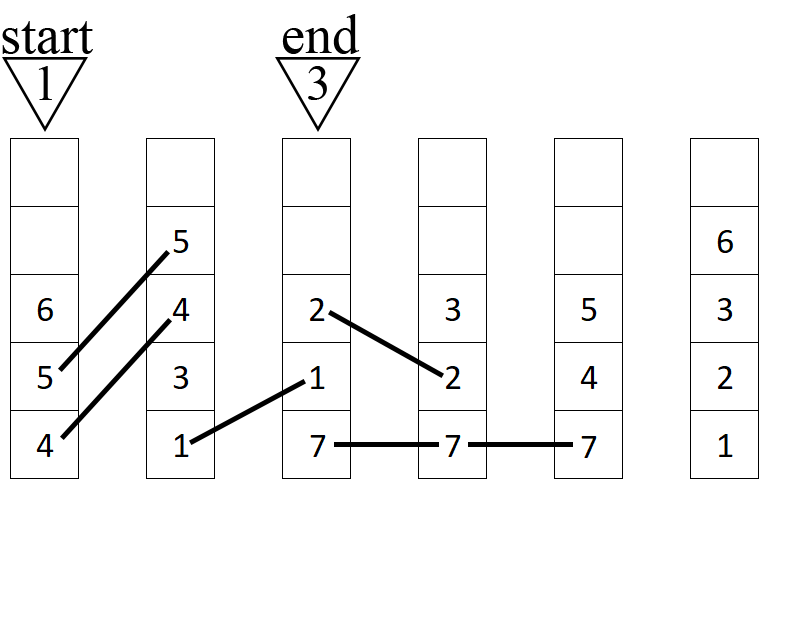}
    \end{minipage} &

    Let's try to build all possible pipes with end at $ end = 3 $ and beginning at $ start = 1 $. There are no such pipes, since there are no tools that are needed both for the job scheduled at the time of $ 1 $ and for the job planned at the time of $ 3 $, while the time $ 2 $ is not needed for the interim memorial, i.e. $(T_1 \cap T_3) \backslash T_2 = \varnothing$.
    \\
	\hline
	
    \begin{minipage}{0.4\textwidth}
      \includegraphics[width=\linewidth]{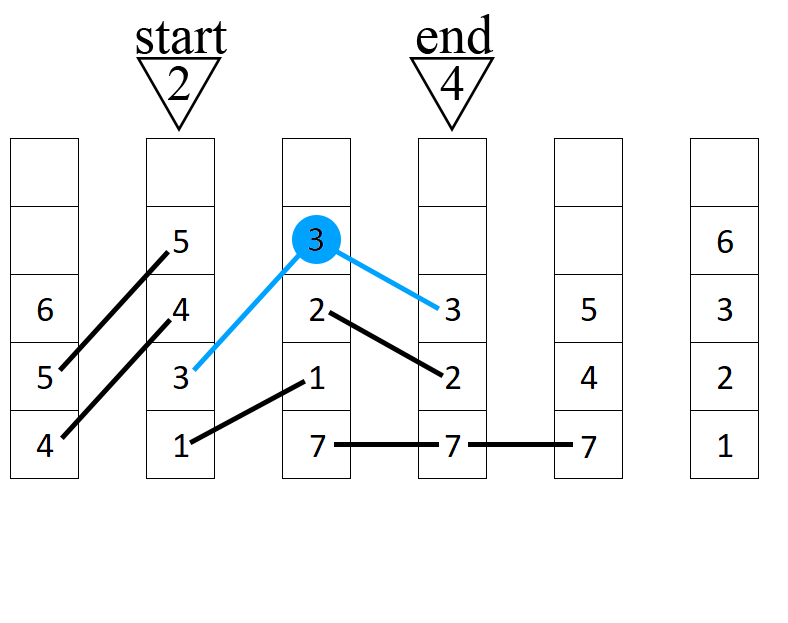}
    \end{minipage} &
Let's try to build all possible pipes with end at $ end = 4 $ and beginning at $ start = 2 $. $ (T_4 \cap T_2) \backslash T_3 = \{3 \} $. At the instant $ 3 $ there are two empty slots, one of which we will occupy with the tool $ 3 $, thereby constructing a pipe $ \pi_{2,4}^{3} $ with the beginning at the instant $ 2 $, the end at the moment $ 4 $ and tool $ 3 $. Note that pipes with an end at $ end = 3 $ can no longer be built, since each tool in this instant is already the end of a pipe and any expression of the form $ (T_i \cap T_4) \backslash (T_3 \cup T_2 \cup \dots) = \varnothing $ since $ T_4 \subseteq T_3 \cup T_2 $. Then the iteration where $ start = 1 $ can be skipped, which is reflected in the \hyperref[theorem:GPCA2]{MPCA-bitwise} in line $ 9 $, where the variable $ | end \_tools | $ will be equal to zero.
    \\
	\hline

  \end{tabular}
\label{table:MPCA_example_1}
\end{table}

\begin{table}[H]

  \caption{
  Example of $MPCA$ execution.
  }

  \begin{tabular}{ | c | m{9cm} | }
  
      \hline

    \begin{minipage}{0.4\textwidth}
      \includegraphics[width=\linewidth]{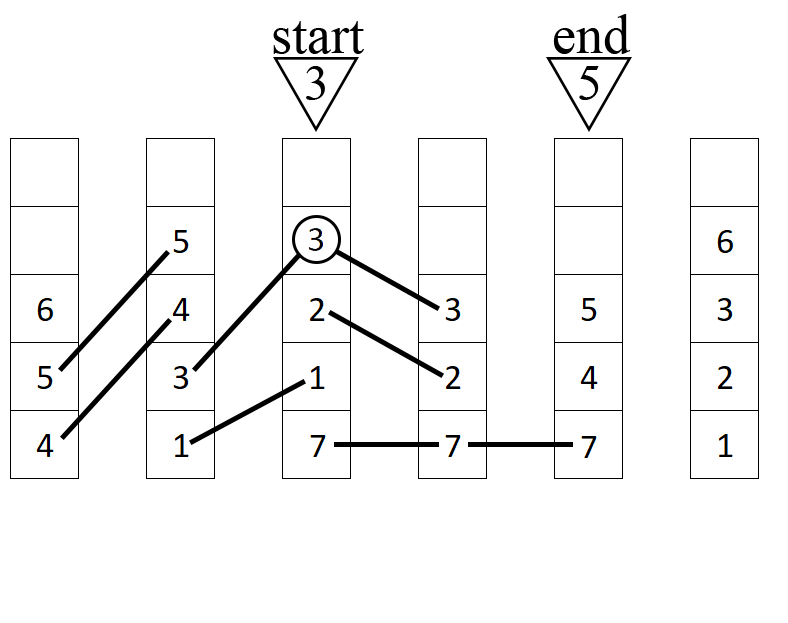}
    \end{minipage} &
Let's try to build all possible pipes with end at $ end = 5 $ and beginning at $ start = 3 $.
There are no such pipes, since there are no tools that are needed both for the job planned at the time of $ 3 $ and for the job planned at the time of $ 5 $, while the time $ 4 $ is not needed for the intervening instant, that is, $ (T_5 \cap T_3) \backslash T_4 = \varnothing $. \\
    
	\hline
	
  	\hline
    \begin{minipage}{0.4\textwidth}
      \includegraphics[width=\linewidth]{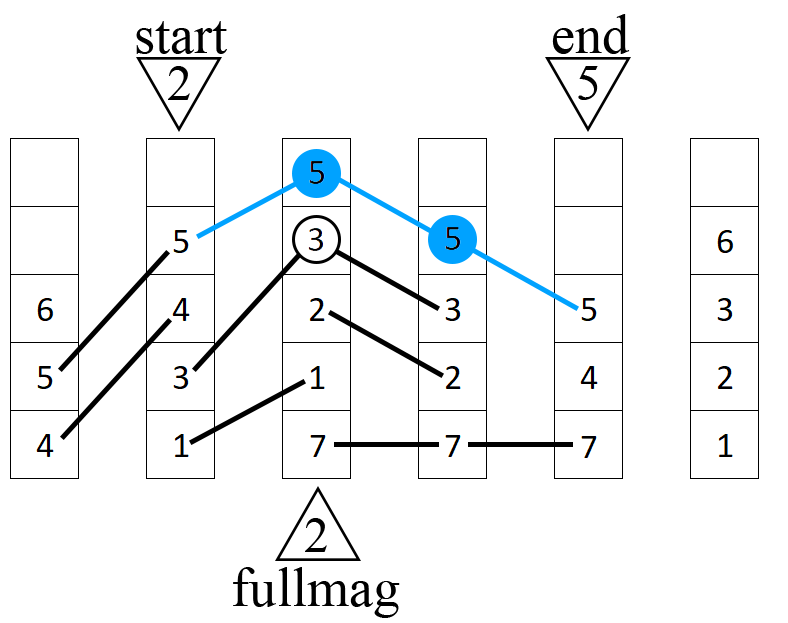}
    \end{minipage} &
Let's try to build all possible pipes with end at $ end = 5 $ and beginning at $ start = 2 $. $ (T_2 \cap T_5) \backslash (T_3 \cup T_4) = \{4,5 \} $. Then the pipes $ \pi_{2,5}^{4} $ and $ \pi_{2,5}^{5} $ claim to be constructed. At instant $ 4 $ there are two empty slots, At instant $ 3 $ there is one empty slot, then only one pipe can be built instead of two. Let's take one empty slot at these instants $ 3.4 $, thereby constructing a pipe $ \pi_{2,5}^{5} $ with the beginning at the instant $ 2 $, the end at the instant $ 5 $ and the tool $ 5 $. Note that at the instant $ 3 $ the magazine is completely full, so no pipes passing through the instant $ 2 $ will already be built. Then the iteration $ start = 1 $ can be skipped, which is \hyperref[theorem:GPCA2]{MPCA-bitwise} in the algorithm in line $ 9 $, when $ fullmag = 2> 1 = start $.
    \\
  
    \hline

    \begin{minipage}{0.4\textwidth}
      \includegraphics[width=\linewidth]{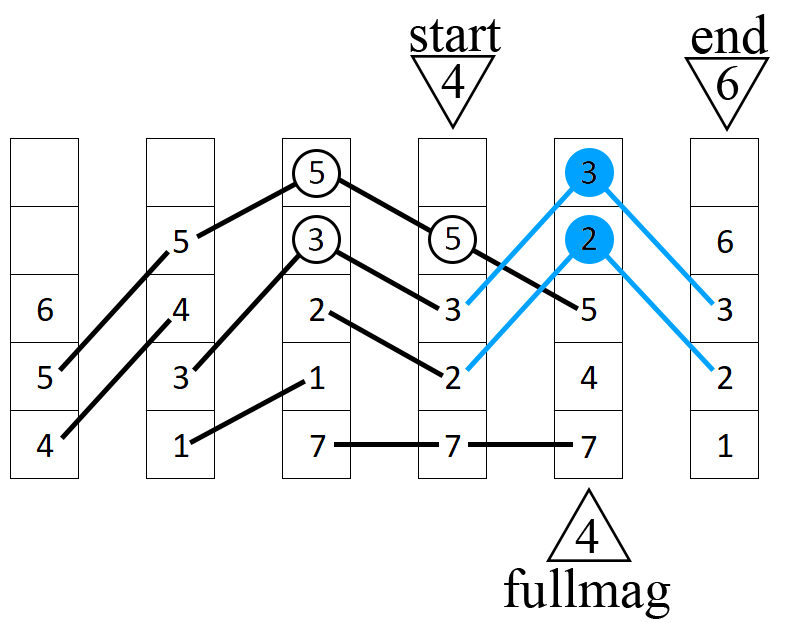}
    \end{minipage} &
Let's try to build all possible pipes with end at $ end = 6 $ and beginning at $ start = 4 $. $ (T_4 \cap T_6) \backslash T_5 = \{2,3 \} $. Then the pipes $ \pi_{4,6}^{2} $ and $ \pi_{4,6,3} $ claim to be constructed. At the instant $ 5 $ there are two empty slots, we will occupy these slots with tools $ 2,3 $, thereby constructing two pipes $ \pi_{4,6,2} $ and $ \pi_{4,6,3} $. Note that at the moment of time $ 5 $ the magazine is completely full, so no pipes passing through the instant $ 5 $ can no longer be built. Then iterations $ start = 3,2,1 $ can be skipped, which is \hyperref[theorem:GPCA2]{MPCA-bitwise} in the algorithm in line $ 9 $. MPCA has finished job. There is $10$ pipes than number of switches = $\sum_{i=1}^{n}|T_i| -C - 10 = 20 - 5 - 10 = 5.$
    \\
	\hline
	
  \end{tabular}
\label{table:MPCA_example_2}
\end{table}

\hyperref[algorithm:MPCA-bitwise]{\textit{MPCA-bitwise}} is an efficient implementation of $ MCPA $, where the sets $ T_i $ and $ M_i $ are encoded by $ 64 $ -bit vectors, which will allow performing set operations in $ \lceil \frac {m} {64} \rceil $ operations.
Let's analyze the complexity. Let's create a sparce table called $possib\_tools$. Time to create it $ O (n log (n)) $, time to call it $ O (1) $. Let $ l $ be the maximum pipe length. The variable $ end\_tools $ symbolizes the set of $ tool $ tools such that the pipe $ \pi_{start, end}^{tool} $ can be built, where $start \geq fullmag$. From this variable with each new $start$ all pipes $\pi_{start,end}^{tool}$ are removed, than when at the moment when $start=end-l$ variable $end\_tools$ will become equal to empty set and loop from string $8$ will finish job because of  $break$ from line $10$.  The loop on line $ 3 $ does $ n-1 $ iterations, the loop on line $ 8 $ does no more than $ l $ iterations. The $ 11 $ line runs in $\lceil \frac{m}{64} \rceil$ operations. Each execution of the line $ 21 $ means loading one slot of the magazine, of which there are only $ C n $, thus $ 21 $ will be called at most $ C n $. Then the time complexity $O(l \lceil \frac{m}{64} \rceil n + n log(n)) = O(l m n + n log(n))$.

\begin{algorithm}[H]
\SetAlgoLined
\DontPrintSemicolon
\SetKw{KwGoTo}{go to}
	$pipes\_count:= 0$\;
	$pissib\_tools[i][j]:=\bigcup_{k=i}^{j} T_k$ for all $i<j$\;
	\For{$end=2,\dots,n$}{
		\If{$empty[end-1] = 0$}{
			$fullmag:=end-1$\;
		}
		$end\_tools:= T_{end} \cap pissib\_tools[fullmag][end-1]$\;	
		\For{$start=end-1,end-3,\dots,1$}{
			\If{$fullmag > start $ $ or $ $ |end\_tools|=0$}{
				\textbf{break}\;
			}
			$candidates:= T_{start} \cap end\_tools$\;
			\If{$|candidates|>0$}{
				$end\_tools:= end\_tools \backslash candidates$\;
				$new\_pipes\_count:= |candidates|$\;
				\If{$new\_pipes\_count > 0$}{
					\For{$j = start+1,start+2,\dots,end-1$}{
						\If{$empty[j] \leq new\_pipes\_count$}{
							$new\_pipes\_count:= empty[j]$\;
							$fullmag:= j$\;
						}
						$empty[j]:= empty[j] - new\_pipes\_count$\;
					}
					$pipes\_count:=pipes\_count + new\_pipes\_count$\;
				}
			}    		
    	}    	
    	
    }
    \Return $pipes\_count$\;

  \caption{$MPCA-bitwise$ }
\end{algorithm}\label{algorithm:GPCA2}

\begin{center}
\textbf{EXPERIMENTAL RESULTS}
\end{center}
To compare the speed of computing the SSP objective function we present our computational study of $ KTNS $, $ MPCA $, $ KTNS-bitwise $, $ MPCA-bitwise $ algorithms as an intermediate implementations from $ KTNS $ to $ MPCA-bitwise $. 

All computations were performed on an Intel\textsuperscript{®} Core\textsuperscript{TM} i$5$ CPU $2.60$ GHz computer with 4 GB or RAM.
$ MPCA-bitwise $, $ KTNS-bitwise $, $ MPCA $ are implemented in $C++$. The implementation of $KTNS$ was taken from the repository published by Mecler et al.{[\ref{Mecler}]} Both MPCA and KTNS was compilated with g++ version 10.3.0 using $-O3$ flag.

\begin{figure}[H]
\begin{minipage}[t]{.5\linewidth}
  \includegraphics[width=\linewidth]{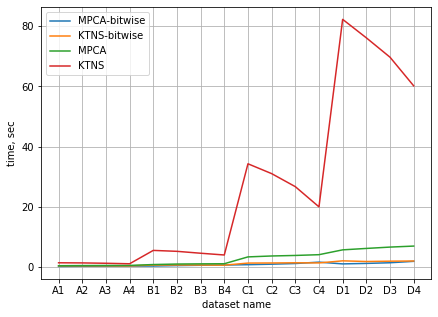}\label{figure:KTNS_vs_MPCA_1}
  \caption{Comparison of KTNS[\ref{Mecler}], MPCA, KTNS-bitwise, MPCA-bitwise for Catanzaro et. al. datasets.}
\end{minipage}
\begin{minipage}[t]{.5\linewidth}
  \includegraphics[width=\linewidth]{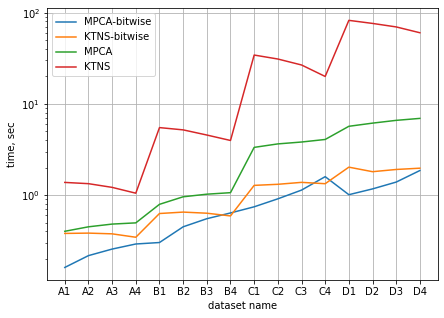}\label{figure:KTNS_vs_MPCA_2}
  \caption{Comparison of KTNS[\ref{Mecler}], MPCA, KTNS-bitwise, MPCA-bitwise for Catanzaro et. al. datasets in logarithmic scale.}
\end{minipage}
\end{figure}

$10^5$ random job sequences were generated for each Catanzaro et al.[\ref{Catanzaro}] dataset, each dataset contains $10$ instances, so each row in \hyperref[table:KTNS_vs_MPCA]{Table 1} shows the processing time of $10^6$ problems by KTNS and proposed MPCA. Results of computational experiments is given in \hyperref[table:MPCA_vs_KTNS]{Table 7},  \hyperref[figure:KTNS_vs_MPCA_1]{Figure 1,2}. You can see, as mentioned in the previous section, $ KTNS $ is accelerating with $ C $ growth, while $ MPCA $ is decelerating, which can be seen in the graph as non-monotonicity of $ KTNS $. On average, MPCA is 6 times faster than KTNS on type A datasets, 11 times on type B datasets, 28 times on type C datasets, 59 times on type D datasets.

\begin{table}[H]
\caption{Comparison of KTNS[\ref{Mecler}], MPCA, KTNS-bitwise, MPCA-bitwise,  for Catanzaro et. al. datasets.}

\begin{tabular}{l|c|c|c|c|c|c|c}

dataset &   n &   m &   C &    KTNS[\ref{Mecler}] &   MPCA &  KTNS-bitwise &  MPCA-bitwise \\
    \hline
          A1 &  10 &  10 &   4 &   1.377 &  0.401 &         0.380 &         0.161 \\
          A2 &  10 &  10 &   5 &   1.334 &  0.449 &         0.383 &         0.217 \\
          A3 &  10 &  10 &   6 &   1.215 &  0.481 &         0.376 &         0.256 \\
          A4 &  10 &  10 &   7 &   1.049 &  0.496 &         0.345 &         0.291 \\
    \hline
          B1 &  15 &  20 &   6 &   5.493 &  0.792 &         0.628 &         0.302 \\
          B2 &  15 &  20 &   8 &   5.187 &  0.958 &         0.651 &         0.449 \\
          B3 &  15 &  20 &  10 &   4.554 &  1.023 &         0.633 &         0.551 \\
          B4 &  15 &  20 &  12 &   3.969 &  1.063 &         0.591 &         0.637 \\
    \hline
          C1 &  30 &  40 &  15 &  34.291 &  3.335 &         1.278 &         0.744 \\
          C2 &  30 &  40 &  17 &  31.021 &  3.640 &         1.313 &         0.910 \\
          C3 &  30 &  40 &  20 &  26.690 &  3.817 &         1.379 &         1.134 \\
          C4 &  30 &  40 &  25 &  19.972 &  4.071 &         1.333 &         1.592 \\
    \hline
          D1 &  40 &  60 &  20 &  82.308 &  5.677 &         2.024 &         1.010 \\
          D2 &  40 &  60 &  22 &  76.167 &  6.149 &         1.805 &         1.169 \\
          D3 &  40 &  60 &  25 &  69.728 &  6.596 &         1.909 &         1.389 \\
          D4 &  40 &  60 &  30 &  60.206 &  6.931 &         1.972 &         1.859 \\
    \hline
\end{tabular}

\label{table:MPCA_vs_KTNS}
\end{table}

\newpage

\begin{center}
\textbf{Summary and Future Research Directions}
\end{center}
Our $ MPCA-bitwise $ algorithm speeds up $ 59 $ times on average compared to $ KTNS $ for large-scale datasets type D [\ref{Catanzaro}]. In further research we aim to obtain a more accurate time complexity of the algorithm and test it on larger SSP benchmark instances. We also intend to measure the effect of incorporating $MPCA-bitwise$ into exact and approximate algorithms for solving SSP. Another goal is to investigate how the pipe characteristics change when increasing the capacity of the magazine $C$, number of jobs $n$, number of tools $m$, and how this correlates with the optimality of the sequence.

\begin{center}
\textbf{PROOF OF THEOREM 1}
\end{center}
Let $S = (1, \dots n)$  be a sequence of jobs, $M \in \boldsymbol {M} (S)$ - sequence of magazine states, such that at instant $i$ job $i$ can be performed i.e. $T_i \subseteq M_i$, where $i = 1, \dots , n$.\\
Let $G_M=(V,A)$ denote graph where $V = \{v_{i}^{tool}: tool\in M_{i}\}$ and $A = \{(v_{i}^{tool},v_{i+1}^{tool}): v_{i}^{tool},v_{i+1}^{tool} \in V\}$ i.e. arc exists iff tool is planned at instant $i$ and at next instant $i+1$. $V$ shows the content of each slot of the magazine at each moment of time, since $v_{i}^{t} \in V$ iff a tool $t$ is contained in magazine at instant $i$. $(v_{i}^{t},v_{i+1}^{t}) \in A$ iff no switch of a tool $t$ at instant $i+1$ is needed. Therefore number of tool switches in $M$ is equal to $C (n-1) - |A|$ e.i. the number of all possible places where switch might be needed minus number of places where switch are not needed.\\

\begin{figure}[H]
\begin{center}
      \includegraphics[width=0.5\linewidth]{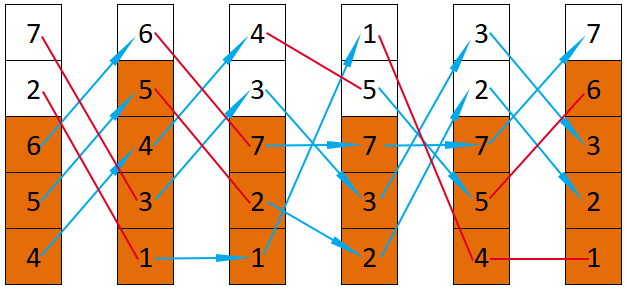}
  \caption{Example of $G_M$. Blue arcs are arcs of $G_M$, red lines symbolize switches.}
\end{center}
\end{figure}

\noindent Let's call a needed vertex a vertex $v_{i}^{t}$ such that $v_{i}^{t} \in T_i$. Let's call a useless vertex a vertex $v_{i}^{t}$ such that $v_{i}^{t} \notin T_i$. \\
$\mathscr{T}(M):=\{\pi_{s,e}^{t}=(v_{s}^{t}, \dots , v_{e}^{t}): t \in T_{s} \cap T_{e}, t \notin \bigcup_{i=s+1}^{e-1} T_{i}, t \in \bigcap_{i=s+1}^{e-1} M_{i} \}$ - the set of all pipes in $M$.\\
Let $\mathscr{H}^{*}(M):= \{P $ - path in  $G_M : P $ contains no common arcs with pipes  $\} $\\
 $\mathscr{H}(M):= \{P \in \mathscr{H}^{*}(M): P$  -  inclusion-wise maximal path $\}$\\
$\mathscr{H}(M) = \mathscr{H}_{0}(M) \sqcup \mathscr{H}_{1}(M) \sqcup \dots \sqcup \mathscr{H}_{n}(M)$, where $\mathscr{H}_{r}(M):=\{P \in \mathscr{H}(M):$ there are exactly $r$ needed vertices $\}$\\

\noindent  \textbf{Lemma 1}: $\mathscr{H}(M)=\mathscr{H}_{0}(M) \sqcup \mathscr{H}_{1}(M)$\\
\textbf{Proof}: \\
Let $P \in \mathscr{H}_{r}(M)$, $r \geq 2$, $P=(v_{1}^{t}, \dots , v_{p}^{t})$ then let's consider two needed vertices $v_{i}^{t}$,$v_{j}^{t}$ such that $v_{i+1}^{t}, \dots , v_{j-1}^{t}$ are useless, then $P$ is a pipe $\pi_{i,j}^{t}$, therefore $P$ contains a pipe as subpath, therefore $P \notin \mathscr{H}^{*}(M)$, therefore $P \notin \mathscr{H}(M)$ $\Box$. \\

\noindent  \textbf{Lemma 1.2.}
$\forall P \neq  P' \in \mathscr{H}(M)$ $A(P) \cap A(P')=\varnothing$\\
\textbf{Proof}:
Suppose $P \neq P'$ and $\exists $ $ a \in A(P) \cap A(P')$, then $P$, $P'$ belong to one connected subgraph of $G_M$. Since the arcs are only between adjacent instants and at one instant the tool cannot be present twice, then any connected subgraph of $G_M$ is a path. $P$, $P' \in \mathscr{H}(M)$ then they are inclusion-wise maximal paths, then if they have a common arc $a$, then $P=P'$, which leads to a contradiction $\Box$.\\

\noindent  \textbf{Lemma 1.3.}
$\forall \pi \neq  \pi' \in \mathscr{T}(M)$ $A(\pi) \cap A(\pi')=\varnothing$\\
\textbf{Proof}:
Suppose $\pi_{s,e}^{t} \neq \pi_{s',e'}^{t'}$ and $\exists a \in A(\pi_{s,e}^{t}) \cap A(\pi_{s',e'}^{t'})$, then $\pi_{s,e}^{t}$, $\pi_{s',e'}^{t'}$ belong to one connected subgraph of $G_M$. Since $\pi_{s,e}^{t}$ and $\pi_{s',e'}^{t'}$ have common arc $a$, than $t=t'$ (Since arcs only connect vertices with the same tool).
Since $\pi_{s,e}^{t}$ is a pipe, then $v_{s+1}^{t}, \dots v_{e-1}^{t}$ are useless vertices, then since $v_{s'}^{t}$ is needed vertex then $v_{s'}^{t}=v_{s}^{t}$ or $v_{s'}^{t} = v_{e}^{t}$. Since $\pi_{s,e}^{t}$ and $\pi_{s',e'}^{t'}$ have common arc $a$, then $\pi_{s',e'}^{t'}$ starts no later than $e-1$, with implies $s=s'$. Similar reasoning leads to $e=e'$.
$\pi_{s,e}^{t}$ and $\pi_{s',e'}^{t'}$ have the same start $s$, the same end $e$, the same tool $t$ and they belong to one connected subgraph of $G_M$ which is path, then  $\pi_{s,e}^{t} = \pi_{s',e'}^{t'}$, which leads to a contradiction $\Box$.\\

\noindent  \textbf{Lemma 1.4.} $switches(M)=\sum_{j=1}^{n}|J_j| -C - |\mathscr{T}(M)|+ |\mathscr{H}_{0}(M)| $\\
\textbf{Proof}:\\
According to the \textbf{Lemma 1} $A(G_M)=A(\mathscr{T}(M)) \sqcup A(\mathscr{H}(M))=A(\mathscr{T}(M)) \sqcup A(\mathscr{H}_{0}(M)) \sqcup A(\mathscr{H}_{1}(M))$\\
$\lambda:=C \cdot n - \sum_{i=1}^{n}|T_i|$ i.e. total number of useless vertices $G_M$.\\
$\lambda=\lambda_{0} + \lambda_1 + \lambda_{2}$, where $\lambda_0$ is number of useless vertices contained in paths from $\mathscr{H}_{0}(M)$, $\lambda_{1}$ from $\mathscr{H}_{1}(M)$, $\lambda_{2}$ from $\mathscr{T}(M)$. \\
Let $P \in \mathscr{H}_{0}(M)$, then
$|A(P)|=|V(P)|-1
= \#($useless vertices in $P) + \#($needed vertices in  $P) -1 
= \#($useless vertices in $P) + 0 -1 = \#($useless vertices in $P) - 1$\\
Then according to the \textbf{Lemma 2} $|A(\mathscr{H}_{0}(M))|= \sum_{P \in \mathscr{H}_{0}(M)} (\#($useless vertices in $P) - 1) = \lambda_0 - |\mathscr{H}_{0}(M)|$\\
Let $P \in \mathscr{H}_{1}(M)$, then
$|A(P)|=|V(P)|-1
= \#($useless vertices in $P) + \#($needed vertices in  $P) -1 
= \#($useless vertices in $P) + 1 -1 = \#($useless vertices in $P)$\\
Then according to the \textbf{Lemma 2} $|A(\mathscr{H}_{1}(M))| = \sum_{P \in \mathscr{H}_{1}(M)} \#($useless vertices in $P) = \lambda_1$.\\
Let $P \in \mathscr{T}(M)$, then $|A(P)|=|V(P)|-1 =\#($useless vertices in $P) + \#($needed vertices in  $P) -1 
= \#($useless vertices in $P) + 2 -1 = \#($useless vertices in $P) +1$ \\
Then according to the \textbf{Lemma 3} $|A(\mathscr{T}(M))| = \sum_{\pi \in \mathscr{T}(M))} (\#($useless vertices in $P)+1) = \lambda_2 + |\mathscr{T}(M)|$.\\
Then $|A(G_M)| = |A(\mathscr{T}(M)) \sqcup A(\mathscr{H}_{0}(M)) \sqcup A(\mathscr{H}_{1}(M))| = |A(\mathscr{H}_{0}(M))| + |A(\mathscr{H}_{1}(M))| + |A(\mathscr{T}(M))|= \lambda_0 - |\mathscr{H}_{0}(M)| + \lambda_1+ \lambda_2 + |\mathscr{T}(M)|=\lambda+|\mathscr{T}(M)|- |\mathscr{H}_{0}(M)| = C \cdot n - \sum_{j=1}^{n}|T_j| + |\mathscr{T}(M)|- |\mathscr{H}_{0}(M)|$.\\
Since $switches(M) = C \cdot (n-1) - |A(G_M)|$, then
$switches(M) = C \cdot (n-1) - (C \cdot n - \sum_{j=1}^{n}|T_j| + |\mathscr{T}(M)|- |\mathscr{H}_{0}(M)|) =\sum_{j=1}^{n}|T_j| -C - |\mathscr{T}(M)|+ |\mathscr{H}_{0}(M)|$ $\Box$.  \\

\noindent  \textbf{Theorem 1.} $\underset{M \in \boldsymbol{M}(S)}{\mathrm{min}}\{switches(M)\} = \sum_{i=1}^{n}|T_i| -C - \underset{M \in \boldsymbol{M}(S)}{\mathrm{max}}\{|\mathscr{T}(M)|\}.$\\
\textbf{Proof}:\\
Since according to \textbf{Theorem 4.} $switches(M)=\sum_{j=1}^{n}|T_j| -C - |\mathscr{T}(M)|+ |\mathscr{H}_{0}(M)|$ \\
Let's proof that $\exists M \in \underset{M \in \boldsymbol{M}(S)}{\mathrm{argmax}}\{|\mathscr{T}(M)|\}$ such that $|\mathscr{H}_{0}(M)|=0$.\\
Let's consider an arbitrary $M \in \underset{M \in \boldsymbol{M}(S)}{\mathrm{argmax}}\{|\mathscr{T}(M)|\}$. Let $R_i = \{ t: \exists $ path $P \in \mathscr{H}_{0}(M)$ such that $v_{i}^{t}$ in $P \}$. \\
Than let $M_1' = M_1 \backslash R_1, \dots , M_n' = M_n \backslash R_n$ i.e. there is $\sum_{i=1}^{n} |R_i|$ empty slots in sequence of magazine states $M'$.\\
If it is proved that it is always possible to fill in an empty slot at instant $i$ so that the added tool $t \in M_{i-1}'$ or $t \in M_{i+1}'$, then after $\sum_{i=1}^{n} |R_i|$ repetitions of such a procedure all slots will be filled and $\mathscr{H}_{0}(M)$ will be empty.\\
Suppose that there is no $M_i$ with empty slot that $t \notin M_i$ and $t \in M_{i-1} \cup M_{i+1}$, then $M_{i-1} \cup M_{i+1} \subseteq M_i$, then $M_{i+1}$ has empty slot too, than $M_{i} \cup M_{i+2} \subseteq M_{i+1}$, then $M_{i}=M_{i+1}$. Then $M_{1} = M_{2} = \dots = M_{n}$  and $|M_1|<C$, than $|\bigcup_{i=1}^{n} M_i| = |M_i| <C$, but $T_i \subseteq M_i'$, $i=1 , \dots ,n$, then $\bigcup_{i=1}^{n} T_i \subseteq \bigcup_{i=1}^{n} M_i$, then $|\bigcup_{i=1}^{n} T_i|<C$, but $|\bigcup_{i=1}^{n} T_i| = m > C$, which leads to a contradiction.\\
Finally $\exists M^{*} \in \underset{M \in \boldsymbol{M}(S)}{\mathrm{argmax}}\{|\mathscr{T}(M)|\}$ such that $|\mathscr{H}_{0}(M^{*})|=0$, then since $\sum_{j=1}^{n}|T_j|$ and $C$ are constants and $|\mathscr{H}_{0}(M)|$ can be decreased to zero, then  $\underset{M \in \boldsymbol{M}(S)}{\mathrm{min}}\{switches(M)\} = \sum_{i=1}^{n}|T_i| -C - \underset{M \in \boldsymbol{M}(S)}{\mathrm{max}}\{|\mathscr{T}(M)|\}$ $\Box .$\\

\begin{center}
\textbf{PROOF OF THEOREM 2}
\end{center}

Let $ \boldsymbol {L} (S) = \{L = (L_ {S_1}, \dots, L_ {S_n}): \forall i \in J $ $ T_i \subseteq L_i, | L_i | \geq C \} $ is a magazine state sequence in which job can be performed in the order $ S \in \boldsymbol {S} $, where empty slots are allowed (which limits the conditions $ | L_i | \geq C $). Note that $ \boldsymbol {M} (S) \subseteq \boldsymbol {L} (S) $, since empty slots are not allowed in $ \boldsymbol {M} (S) $, i.e. the magazine will pay in full at every moment of time, i.e. $ \forall i \in J $ $ | M_i | = C $. \\

Let $possib\_pipes(S) = \{\pi_{s,e}^{t}: \exists L \in \boldsymbol{L}(S):\pi_{s,e}^{t} \in \mathscr{T}(L)\} = \{\pi_{s,e}^{t}: |T_j|<C, j=s+1, \dots , e-1, t \in (T_{s} $ $\cap$ $ T_{e}) \backslash (\bigcup_{i=s+1}^{e-1} T_{i}) \}$ denote all possible pipes that can be constructed(not at the same time),  where $S$ is a sequence of jobs. \\
Let $can\_construct(L) = \{\pi_{s,e}^{t} \in possib\_pipes(S): t \notin L_j,|L_j|<C, j=s+1, \dots , e-1\}$ denote all pipes that can be constructed in $L$. \\
We denote by $ \boldsymbol {\widehat {L}} (S) $ the set of all possible sequences of the state of the magazine with possible empty slots, empty slots were filled only when constructing pipes. \\
\textbf{Lemma 2.1}: Let $L \in \boldsymbol{\widehat{L}}(S)$ then $can\_construct(L) = \{\pi_{s,e}^{t} \in possib\_pipes(S) \backslash \mathscr{T}(L): \forall j \in \{s+1, \dots , e-1\}  $ $ |L_j|<C \}$.\\
That is if only pipes has been constructed in $L$, than pipe $\pi_{s,e}^{t}$ can be constructed iff $\pi_{s,e}^{t}$ has not been constructed in $L$ and there is enough empty slots to construct $\pi_{s,e}^{t}$.\\
\textbf{Proof}:\\
Based on the definition of $can\_construct( \cdot)$, it suffices to prove that if $L \in \boldsymbol{\widehat{L}}(S)$ and $\pi_{s,e}^{t} \in possib\_pipes(S) \backslash \mathscr{T}(L)$, then $\forall j \in \{s+1, \dots , e-1\}  $ $ t \notin L_j$. Suppose $\exists j \in \{s+1, d\dots , n \}: t \in L_j$, since $\pi_{s,e}^{t}$ is a pipe, then $\forall i \in \{s+1,\dots , e-1 \}  $ $ t \notin T_i$, , then $t \notin T_j$, which implies that vertex $v_{j}^{t}$ is useless, then since only pipes has been constructed $\exists \pi_{s',e'}^{t} \in \mathscr{T}(L): s' < j < e'$. Since $\pi_{s',e'}^{t}$ is a pipe, than $\forall i \in \{s'+1,\dots , e'-1 \}  $ $ t \notin T_i$, and since $\pi_{s,e}^{t}$ is a pipe, than $ $ $ t \in T_s$, than $s \leq s'$. Since $\pi_{s,e}^{t}$ is a pipe, than $\forall i \in \{s+1,\dots , e-1 \}  $ $ t \notin T_i$, and since $\pi_{s',e'}^{t}$ is a pipe, than $ $ $ t \in T_{s'}$ than $s' \leq s$, then $s=s'$ and similarly $e=e'$, which implies that $\pi_{s',e'}^{t}$ and $\pi_{s,e}^{t}$ are the same pipe, then $\pi_{s,e}^{t} \in \mathscr{T}(L)$ and $\pi_{s,e}^{t} \in possib\_pipes(S) \backslash \mathscr{T}(L)$, which leads to a contradiction $\Box$.\\
Further we always assume that $ L \in \boldsymbol {\widehat {L}} (S) $, since we will only talk about constructing and removing pipes and otherwise empty slots will not be filled. \\
\noindent  Let $ \boldsymbol {L} _ {opt} (S) = \{L \in \boldsymbol {\widehat {L}} (S): | \mathscr {T} (L) | = max \{| \mathscr {T} (M) |: M \in \boldsymbol {M} (S) \} \} $, i.e. this is the set $ L \in \boldsymbol {\widehat {L}} (S) $: $ L $ contains the smallest possible number of pipes. \\
Let $ L \notin \boldsymbol {L_ {opt}} (S) $, then $ L $ does not contain the maximum possible number of pipes, i.e. you can remove some (possibly all) constructed pipes, let them be $ r_1 $ pieces, and then build another set of pipes of capacity $ r_2> r_1 $. $L \notin \boldsymbol{L_{opt}}(S)$ then $ \exists K \subseteq possib\_pipes (S): $ if we remove $ K $ from $ L $, then it will be possible to construct $ K '\subseteq pipes (S): | K' |> | K | $ and vice versa , if there is such a pair $ (K, K ') $, then $ L $ contains not the maximum possible number of pipes, i.e. $ L \notin \boldsymbol {L_ {opt}} (S) $. Let us define $ \mathscr {K} (L) $ as the set of all such pairs $ (K, K ') $. And therefore $\mathscr{K}(L) = \varnothing \iff L \notin \boldsymbol{L_{opt}}(S)$.\\
Let $ \mathscr {K} _ {min} (L): = \{(K, K ') \in \mathscr {K} (L): \nexists (\widetilde{K}, \widetilde{K}' ) \in \mathscr {K} (L): \widetilde{K} \subseteq K $ and $ \widetilde{K}' \subseteq K '$ and $ (\widetilde{K} \neq K $ or $ \widetilde{K}' \neq K ') \} $ i.e. there are many such pairs in which there are no deletions and constructions that could be excluded (for example, it is pointless to delete and add the same pipe). \\
Since $ \mathscr {K} _ {min} (L) = \varnothing \iff \mathscr {K} (L) = \varnothing $, then \\
\textbf {Lemma 2.2}: $ L \in \boldsymbol {L_ {opt}} (S) \iff \mathscr {K} _ {min} (L) = \varnothing $. \\

Let $ needed\_instants (\pi_ {start, end}^{tool}) = \{start + 1, start + 2, \dots, end-1 \} $, i.e. the set of all times in which at least one empty slot is needed to build the pipe $ \pi_ {start, end}^{tool} $.
Let $ needed\_instants (K) = \bigcup _ {\pi_ {s, e}^{t} \in K} needed\_instants (\pi_ {s, e}^{t}) $ i.e. the set of all times in which at least one empty slot is needed to construct the pipe $ \pi_ {s, e}^{t} $: $ \pi_ {s, e}^{t} \in K $. \\

\noindent  \textbf{Lemma 2.3}:
Let $ (K, K ') \in \mathscr {K} (L), \pi_ {s, e}^{t} \in K, \tau_ {s', e '}^{t'} \in K '$, then if\\$needed\_instants(\pi_{s,e}^{t}) $ $\cap$ $ needed\_instants(K') \subseteq needed\_instants(\tau_{s,e}^{t})$, then $(K,K') \notin \mathscr{K}_{min}(L)$\\
\textbf{Proof}:\\
Let's denote $\pi = \pi_{s,e}^{t}$, $\tau = \tau_{s',e'}^{t'}$. 
We get $ L^{0} $ by removing $ K $ from $ L $.
We obtain $ L^{\tau} $ by constructing $ \tau $ in $ L^{0} $.
We obtain $ L^{\pi} $ by constructing $ \pi $ in $ L^{0} $.
$\widetilde{K}':=K' \backslash \{ $ $ \pi \}$. Thus, $ \widetilde {K} '$ can be constructed in $ L^{\pi} $, let us prove that $ \widetilde {K}' $ can also be constructed in $ L^{\tau} $.
Since $needed\_instants(\pi) $ $\cap$ $ needed\_instants(K') \subseteq needed\_instants(\tau)$, then $needed\_instants(\pi) $ $\cap$ $ needed\_instants(K') \subseteq needed\_instants(\tau) $ $\cap$ $ needed\_instants(K')$.  Since $\widetilde{K}' $ $\subseteq$ $ K'$, then \\
$needed\_instants(\pi) $ $\cap$ $ needed\_instants(\widetilde{K}') $ $\subseteq$ $ needed\_instants(\tau) $ $\cap$ $ needed\_instants(\widetilde{K}')$, which (according to \textbf{Lemma 2.1}) means that $\widetilde{K}'$ can be constructed in $ L^{\pi} $, then $ \widetilde {K} '$ can also be constructed in $ L^{\tau} $ that is $ \pi $ can be excluded from the set of $ K $ being removed, and $ \tau $ can be excluded from the set of $ K '$ being added. Then $\widetilde{K}:=K \backslash \{ \pi \}$,$\widetilde{K}':=K' \backslash \{ \tau \}$,then $\exists$ $(\widetilde{K},\widetilde{K}') \in \mathscr{K}(L):\widetilde{K} \subset K $ and $\widetilde{K}' \subseteq K' $, therefore $(K,K') \notin \mathscr{K}_{min}(L)$. $\Box$\\

\noindent 
\textbf{Example}: $needed\_instants(\widetilde{K}'):=\{2,3,4\}$, $needed\_instants(\pi):=\{ 3,4,5,6,7 \}$,\\$needed\_instants(\tau):=\{3,4\}$. According to \textbf{Lemma 2.1} since pipes from $ \widetilde {K} '$ will not use slots whose time points are from $ \{5,6,7 \} $, then for the possibility of constructing $ \widetilde {K}' $, it does not matter whether they are filled or not, it is important only is there enough empty slots at the time instants $ needed\_instants (\widetilde {K} ') $, then the "threat" for $ \widetilde {K}' $ from constructing $ \pi $ is $ needed\_instants (\widetilde { K} ') $ $\cap$ $ needed\_instants (\pi) = \{3,4 \} $, the "threat" to $ \widetilde {K}' $ from building $ \tau $ is $ needed\_instants (\widetilde { K} ') $ $\cap$ $needed\_instants (\tau) = \{3,4 \} $ i.e. with respect to $ \widetilde {K} '$ the construction of $ \pi $ is no worse than the construction of $ \tau $, and if so, if $ \pi $ is contained in the set of $ K $ to be removed, $ \tau $ is contained in the set of $ K '$ to be added, then $ \pi, \tau $ can be simultaneously removed from these sets by obtaining the sets $ (\widetilde {K}, \widetilde {K}') \in \mathscr {K} (L) $, then by definition $ \mathscr {K} _ {min} (L) $ it is true that $ \mathscr {K} _ {min} (L) $ does not contain $ (K, K ') $. \\

\noindent \textbf{Theorem 2}: 
$|\mathscr{T}(MPCA(S))| = max\{ |\mathscr{T}(M)|: M \in \boldsymbol{M}(S) \}$\\
\noindent \textbf{Proof}:\\
According to \textbf {Lemma 2.2}, it will suffice to prove that
$\mathscr{K}_{min}(L)=\varnothing$, where $L=MPCA(S)$\\
Suppose the opposite, i.e. $\exists (K,K') \in \mathscr{K}_{min}(L)$. Since the \hyperref[algorithm:GPCA1]{Algorithm 1} tries to build all pipes from $ possib\_pipes (S) $ and builds if possible. Then it is impossible to complete one more pipe in $ L $ without deleting one before this, therefore the set of pipes to be removed is always not empty. $ K \neq \varnothing $ and since by definition $ | K '|> | K | $, then $ K' \neq \varnothing $. \\
$min\_end(K):=min \{ end : \pi_{start,end}^{tool} \in K \}$\\

\textbf{Case 1}: $e = min\_end(K)=min\_end(K')=e'$.\\
Then let $ \pi_ {s, e}^{t} \in K $, $ \pi_ {s ', e'}^{t '} \in K' $. 2 cases are possible
\begin{enumerate}
\item $s \geq s'$. Then since $e'=e$, then $needed\_instants(\pi_{s,e}^{t}) \subseteq needed\_instants(\pi_{s',e'}^{t'})$. Therefore, according to the Lemma 2.3 $(K,K') \notin \mathscr{K}_{min}(L)$, which leads to a contradiction.
\item $s < s'$. Since $needed\_instants(\pi_{s',e'}^{t'}) \subseteq  needed\_instants(\pi_{s,e}^{t})$, then after deletion $\pi_{s,e}^{t}$  it will be possible to build $\pi_{s',e'}^{t'}$ and therefore it will be possible to move on to $(\widetilde{K},\widetilde{K}'): \widetilde{K}=K \backslash \{ \pi_{s,e}^{t} \}$, $\widetilde{K}'=K' \backslash \{ \pi_{s',e'}^{t'} \}$. It is impossible to endlessly get into \textbf {Case 1.2} because each time the cardinality of $ K $ and $ K '$ decreases by $ 1 $.

\end{enumerate}

\textbf{Case 2}: $e = min\_end(K) < min\_end(K')=e'$.\\
Then let $ \pi_ {s, e}^{t} \in K $, $ \pi_ {s ', e'}^{t '} \in K' $. 2 cases are possible:
\begin{enumerate}
\item $s \geq s'$. Then since $ e <e '$, then $needed\_instants(\pi_{s,e}^{t}) \subseteq needed\_instants(\pi_{s',e'}^{t'})$. Therefore, according to Lemma 2.3 $(K,K') \notin \mathscr{K}_{min}(L)$, which leads to a contradiction.
\item $s < s'$. Let $s''=min\{i: \pi_{i,j}^{k} \in K'\}$ and we will consider the pipe $\pi_{s'',e''}^{t''} \in K'$, then 

since $ s' \geq s''$ either we go to \textbf {Case 2.1}, or $ s <s'' $ and from 
the minimality of $ s''$ and $ e <e''$ it follows that $needed\_instants(\pi_{s,e}^{t}) 
$ $\cap$ $ needed\_instants(K') \subseteq needed\_instants(\pi_{s'',e''}^{t''})$. Therefore, according to L
emma 4 $(K,K') \notin \mathscr{K}_{min}(L)$, which leads to a contradiction.
\end{enumerate}

\textbf{Case 3}: $e = min\_end(K) > min\_end(K')=e'$.\\
Then let $\pi_{s',e'}^{t'} \in K'$, $s=min\{i: \pi_{i,j}^{k} \in K, j=end \}$, $\pi_{s,e}^{t} \in K$.

Emptying the slots at each of the times from the set $ needed\_instants (\pi_ {s', e '}^{t'}) $ $\cap$ $ needed\_instants (K) $ will allow you to build a pipe $ \pi_ {s', e ' }^{t '} $, otherwise even after removing $ K $ it would not be possible to build the pipe $ \pi_ {s', e '}^{t'} $. But note that $ s $ is minimally possible and $ e> e '$, which implies that $ needed\_instants (\pi_ {s', e '}^{t'}) $ $\cap$ $ needed\_instants (K) \subseteq needed\_instants (\pi_ {s, e}^{t}) $, which means that after removing the pipe $ \pi_ {s, e}^{t} $ it will be possible to construct $ \pi_ {s' , e '}^{t'} $. Consequently, the \hyperref[algorithm:GPCA1]{Algorithm 1} constructed the pipe $ \pi_ {s ', e'}^{t '} $ at an iteration earlier than $ \pi_ {s, e}^{t} $, since $ e' <e $, and the \hyperref[algorithm:GPCA1]{Algorithm 1} iterates over the variable $ end $ in ascending order. At the iteration, when the variable $ end $ was equal to $ e '$, the pipe $ \pi_ {s', e '}^{t'} $ should have been built, since at that moment the pipes $ \pi_ {s, e}^{t} $ didn't exist yet, which is the same as it was removed. From which it follows that $ \pi_ {s ', e'}^{t '} $ has already been built, which contradicts the fact that it belongs to $ K' $. Let's move on to  $ (\widetilde {K}, \widetilde {K} '): \widetilde {K} = K \backslash \{\pi_ {s', e '}^{t'} \} $, $ \widetilde {K} '= K' \backslash \{\pi_ {s', e '}^{t'} \} $ and note that after removing $ \widetilde {K} $ it will be possible to build $ \widetilde {K } '$ i.e. $ (\widetilde {K}, \widetilde {K} ') \in \mathscr {K} (L) $, but $ \widetilde {K} \subset K $, $ \widetilde {K}' \subseteq K ' $. Then, by the definition of $ \mathscr {K} _ {min} (L) $, $ (K, K ') \notin \mathscr {K} _ {min} (L) $, which leads to a contradiction.

\textbf {Case 1}, \textbf {Case 2}, \textbf {Case 3} led to a contradiction, hence the assumption about $ \exists (K, K ') \in \mathscr {K} _ {min} (MPCA ( S)) $ is not true, then by Lemma 2.2 we have $ MPCA (S) \in \boldsymbol {L_ {opt}} (S) $, which means that $|\mathscr{T}(MPCA(S))| = max\{ |\mathscr{T}(M)|: M \in \boldsymbol{M}(S) \}$ $\Box$.

\newpage
\begin{center}
\textbf{REFERENCES}
\end{center}

\renewcommand{\refname}{}


\begin{thebibliography}{00}

\bibitem{latexcompanion}\label{TangDenardo} 
C.S. Tang, E.V. Denardo, Models arising from a flexible manufacturing
machine, part I: minimization of the number of tool switches, Oper. Res. 36, \url{http://dx.doi.org/10.1287/opre.36.5.767}.

\bibitem{latexcompanion}\label{daSilva2021}
T. T. da Silva, A. A. Chaves and H. H. Yanasse, A new multicommodity ﬂow model for the job sequencing and tool switching problem. International Journal of Production Research, 59, (2021) 3617–3632, \url{https://doi.org/10.1080/00207543.2020.1748906}

\bibitem{latexcompanion}\label{GoldRom2021}
B. Goldengorin, V. V. Romanuke, Experimental analysis of tardiness in preemptive single machine scheduling. Expert Syst. Appl., 186, (2021) 114947,  \url{ https://doi.org/10.1016/j.eswa.2021.114947}

\bibitem{latexcompanion}\label{AhmadiGoldengorin}
E. Ahmadi, B. Goldengorin, G.A. Süer, H. Mosadegh, A hybrid method of 2-TSP and novel learning-based GA for job sequencing and tool switching problem. Applied Soft Computing, 65,(2018) 214-229.  \url{https://doi.org/10.1016/j.asoc.2017.12.045}

\bibitem{latexcompanion}\label{Catanzaro}
D. Catanzaro, L. Gouveia, M. Labbé, Improved integer linear programming
formulations for the job Sequencing and tool Switching Problem, Eur. J. Oper.
Res. 244 (2015) 766–777, \url{http://dx.doi.org/10.1016/j.ejor.2015.02.018}.

\bibitem{latexcompanion}\label{Calmels}
Calmels, D., 2019. The job sequencing and tool switching problem: State-of-the-art literature review, classification, and trends. International Journal of Production Research 57, 5005–5025
\url{https://doi.org/10.1080/00207543.2018.1505057}

\bibitem{latexcompanion}\label{Mecler}
J. Mecler, A. Subramanian, T. Vidal  A simple and effective hybrid genetic search for the job sequencing and tool switching problem  Computers and Operations Research. 2021. 127, 105153
\url{https://doi.org/10.1016/j.cor.2020.105153}

\bibitem{latexcompanion}\label{Johnson}
D.  Johnson A theoretician's guide to the experimental analysis of algorithms, DIMACS Series in Discrete Mathematics and Theoretical Computer Science. 2002. V. 59
\url{https://doi.org/10.1090/dimacs/059}

\bibitem{latexcompanion}\label{Ghiani}
G. Ghiani, A. Grieco, E. Guerriero, 2010. Solving the job sequencing and tool switching problem as a nonlinear least cost hamiltonian cycle problem. Netjobs 55, 379–385. \url{https://doi.org/10.1002/net.20341}

\bibitem{latexcompanion}\label{Crama}
Y. Crama, A.W.J. Kolen, A.G. Oerlemans, F.C.R. Spieksma, Minimizing the number of tool switches on a flexible machine, Int. J. Flex. Manuf. Syst. 6 (1994) 33–54, \url{http://dx.doi.org/10.1007/BF01324874}



\end{thebibliography}
\end{document}